\documentclass[11pt]{amsart}
\usepackage{amsmath, latexsym, amssymb}
\input psfig.tex

 \newcommand{\DRLambda}{\Tilde{\Tilde{\Lambda}}} 
 \newcommand{\bA}{{\mathbf{ A}}} \newcommand{\bB}{{\mathbf{ B}}} 
 \newcommand{\bb}{{\mathbf{ b}}} \newcommand{\bC}{{\mathbf{ C}}} 
 \newcommand{\bD}{{\mathbf{ D}}}  
 \newcommand{\be}{{\mathbf{ e}}} \newcommand{\bH}{{\mathbf{ H}}} 
 \newcommand{\bF}{{\mathbf{ F}}} \newcommand{\bG}{{\mathbf{ G}}} 
 \newcommand{\bI}{{\mathbf{ I}}} \newcommand{\bJ}{{\mathbf{ J}}} 
 \newcommand{\bK}{{\mathbf{ K}}}  
 \newcommand{\bM}{{\mathbf{ M}}} \newcommand{\bL}{{\mathbf{ L}}} 
 \newcommand{\bN}{{\mathbf{ N}}} \newcommand{\bQ}{{\mathbf{ Q}}} 
 \newcommand{\bq}{{\mathbf{ q}}} \newcommand{\bP}{{\mathbf{ P}}} 
 \newcommand{\bp}{{\mathbf{ p}}}  
 \newcommand{\bR}{{\mathbf{ R}}}  
 \newcommand{\bS}{{\mathbf{ S}}}  
 \newcommand{\bT}{{\mathbf{ T}}}  
 \newcommand{\bU}{{\mathbf{ U}}} \newcommand{\bV}{{\mathbf{ V}}} 
  \newcommand{\bx}{{\mathbf{ x}}} 
 \newcommand{\bY}{{\mathbf{ Y}}} \newcommand{\by}{{\mathbf{ y}}} 
 \newcommand{\bZ}{{\mathbf{ Z}}} \newcommand{\bz}{{\mathbf{ z}}} 
 \newcommand{\bu}{{\mathbf{ u}}} \newcommand{\bv}{{\mathbf{ v}}} 
 \newcommand{\bW}{{\mathbf{ W}}}  
    
 \numberwithin{equation}{section} 
 \newtheorem{theorem}{Theorem}[section] 
  
 \newtheorem{lemma}[theorem]{Lemma} \title{ Variations of Ritz and 
 Lehmann Bounds} \author{Christopher Beattie} \address{Department of 
 Mathematics, Virginia Polytechnic Institute and State University, 
 Blacksburg, VA 24061 USA }

\begin{document}
\begin{abstract} Eigenvalue estimates that are optimal in some sense have 
self-evident appeal and leave estimators with a sense of virtue and 
economy. So, it is natural that ongoing searches for effective strategies 
for difficult tasks such as estimating matrix eigenvalues that are situated 
well into the interior of the spectrum  revisit from 
time to time methods that are known to yield optimal bounds.  This 
article reviews a variety of results related to obtaining optimal bounds to 
matrix eigenvalues --- some results are well-known; others are less 
known; and a few are new.  We focus especially on Ritz and harmonic 
Ritz values, and right- and left-definite variants of Lehmann's 
method.
\end{abstract}
\maketitle

\section{Ritz and Related Values}
 Let $\bK$ and $\bM$ be $n \times n$ real symmetric  positive definite
	matrices and consider the eigenvalue problem
\begin{equation}
	\bK\bx = \lambda \, \bM \bx
	\label{eigprob}
\end{equation}
Label the eigenvalues from the edges toward the center (following  \cite{Parlett})
as  
$$
\lambda_1 \leq \lambda_2 \leq \lambda_{3}\leq \dots
\leq \lambda_{-3} \leq \lambda_{-2} \leq \lambda_{-1}
$$ 
with labeling inherited by the associated eigenvectors:
$\bx_{1},\ \bx_{2},\ \dots,\ \bx_{-2},\ \bx_{-1}$.

Solutions to (\ref{eigprob}) are evidently eigenvalue/eigenvector pairs  of
the matrix $\bM^{-1}\bK$, which is non-symmetric on the face of it.  
However, $\bM^{-1}\bK$ is self-adjoint with respect to the both the
$\bM$-inner product, $\bx^{t}\bM\bx$, and the $\bK$-inner product, 
$\bx^{t}\bK\bx$.  Denote by $\bx^{\mathsf{m}}$ the $\bM$-adjoint of a vector $\bx$, 
$\bx^{\mathsf{m}}=\bx^{t}\bM$, and by $\bx^{\mathsf{k}}$ 
the $\bK$-adjoint, $\bx^{\mathsf{k}}=\bx^{t}\bK$. ``Self-adjointness'' of 
$\bM^{-1}\bK$ amounts to the assertion that for all $\bx$ and $\by$,
$\bx^{\mathsf{m}}(\bM^{-1}\bK\by)=(\bM^{-1}\bK\bx)^{\mathsf{m}}\by$ and 
$\bx^{\mathsf{k}}(\bM^{-1}\bK\by)=(\bM^{-1}\bK\bx)^{\mathsf{k}}\by$.  
Self-adjointness with respect to the $\bM$- and $\bK$-inner products 
implies that the matrix representation of $\bM^{-1}\bK$ with respect to 
any $\bM$-orthogonal or $\bK$-orthogonal basis will be \emph{symmetric}.

 For a given subspace ${\mathcal P}$ of dimension $m < n$, the 
Rayleigh-Ritz method proceeds by selecting a basis for  ${\mathcal P}$, 
say constituting the columns of a matrix 
$\bP\in {\mathbb  R}^{n \times m}$, and then considering the (smaller) 
eigenvalue problem
\begin{equation}
	\bP^{t}\bK\bP\by = \Lambda \ \bP^{t}\bM\bP \by.
	\label{Ritzeig}
\end{equation}
This will yield $m$ eigenvalues (called \emph{Ritz values}) labeled 
similarly to $\{\lambda_{i}\}$ as
$$
\Lambda_1 \leq \Lambda_2 \leq \Lambda_{3}  \leq \dots \leq \Lambda_{-3} 
\leq \Lambda_{-2} \leq \Lambda_{-1}
$$ 
with corresponding eigenvectors 
$\by_{1},\  \by_{2},\ \dots \ \by_{-2},\ \by_{-1}$.
  Vectors in ${\mathcal P}$ given as $\bu_{k}=\bP\by_{ k}$ are 
  \emph{Ritz vectors} associated with the Ritz values $\Lambda_{ k}$.  
  Since $\mathbb{R}^{\mathsf{m}}=\mbox{span}(\mathbf{y}_{1},\ 
  \mathbf{y}_{2},\ \dots,\ \by_{-2},\ \by_{-1}),$ the full set of Ritz 
  vectors evidently forms a basis for ${\mathcal P}$, which  is both
  $\bK$-orthogonal and $\bM$-orthogonal and may be presumed to be 
  $\bM$-normalized without loss of generality: 
  $\bu_{i}^{\mathsf{k}}\bu_{j}=0$ and $\bu_{i}^{\mathsf{m}}\bu_{j}=0$ 
  for $i\neq j$, and $\bu_{i}^{\mathsf{m}}\bu_{i}=1$.

\emph{Harmonic Ritz} values \cite{PaiParVor} result from applying the Rayleigh-Ritz method 
to the eigenvalue problem 
\begin{equation}
\bK	\bM^{-1} \bK\bx = \lambda \,\bK\bx,
	\label{eigprob2}
\end{equation}
which is equivalent to (\ref{eigprob}) --- it has the same eigenvalues and 
eigenvectors.    If we
use the same subspace ${\mathcal P}$, the harmonic Ritz values
are then the eigenvalues of the $m\times m$ problem
\begin{equation}
	\bP^{t}\bK\bM^{-1}\bK\bP\by = \tilde{\Lambda} \, \bP^{t}\bK\bP \by,
	\label{HarmRitzeig}
\end{equation}
yielding  
$$
\tilde{\Lambda}_1 \leq \tilde{\Lambda}_2 \leq \tilde{\Lambda}_{3}\leq  
\dots \leq \tilde{\Lambda}_{-3} \leq \tilde{\Lambda}_{-2} 
\leq \tilde{\Lambda}_{-1}.
$$ Just as Ritz values are weighted means of the 
eigenvalues of the matrix, harmonic Ritz values are harmonic means of the 
eigenvalues of the matrix. 

Quantities which will be introduced here (for lack of a better name) as
\emph{dual harmonic Ritz} values result from applying the Rayleigh-Ritz method 
to the eigenvalue problem 
\begin{equation}
\bM\bx = \lambda \,\bM\bK^{-1}\bM\bx,
\label{eigprob3}
\end{equation}
which is also equivalent to (\ref{eigprob}), in the sense of having the same eigenvalues and 
eigenvectors.  If we
use the same approximating subspace ${\mathcal P}$, the dual harmonic Ritz values
are  the eigenvalues of the $m\times m$ problem
\begin{equation}
\bP^{t}\bM\bP\by = \DRLambda  \, \bP^{t}\bM\bK^{-1}\bM\bP \by,
	\label{DualHarmRitzeig}
\end{equation}
yielding  
$$
\DRLambda_1 \leq \DRLambda_2 \leq \DRLambda_{3} \leq \dots \leq 
\DRLambda_{-3} \leq \DRLambda_{-2} \leq \DRLambda_{-1}.  $$
Dual harmonic Ritz values are also harmonic means of the matrix 
eigenvalues, however with a different weighting than for harmonic Ritz 
values.
 \begin{figure}
  \centerline{\psfig{figure=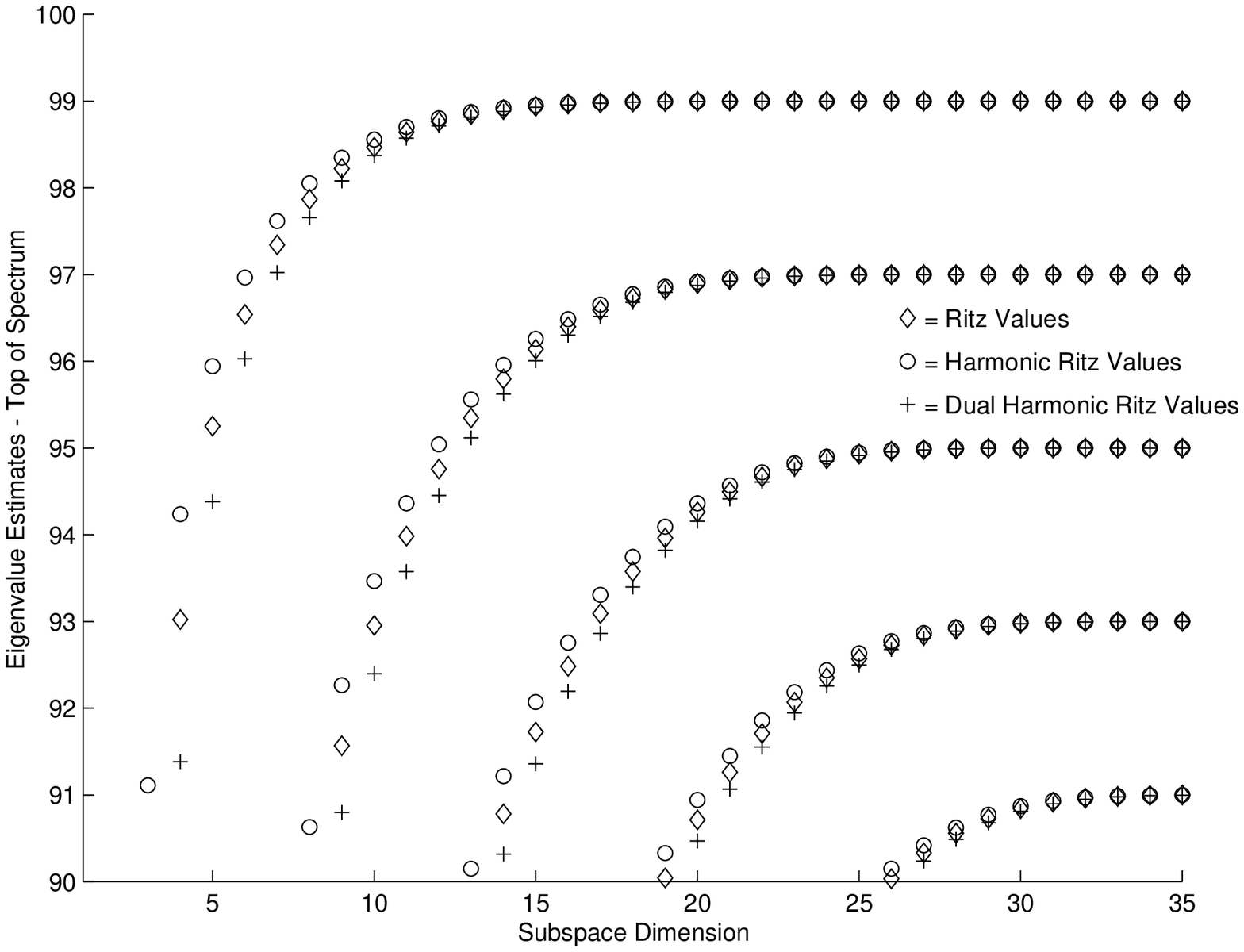,height=3in}}
\centerline{\psfig{figure=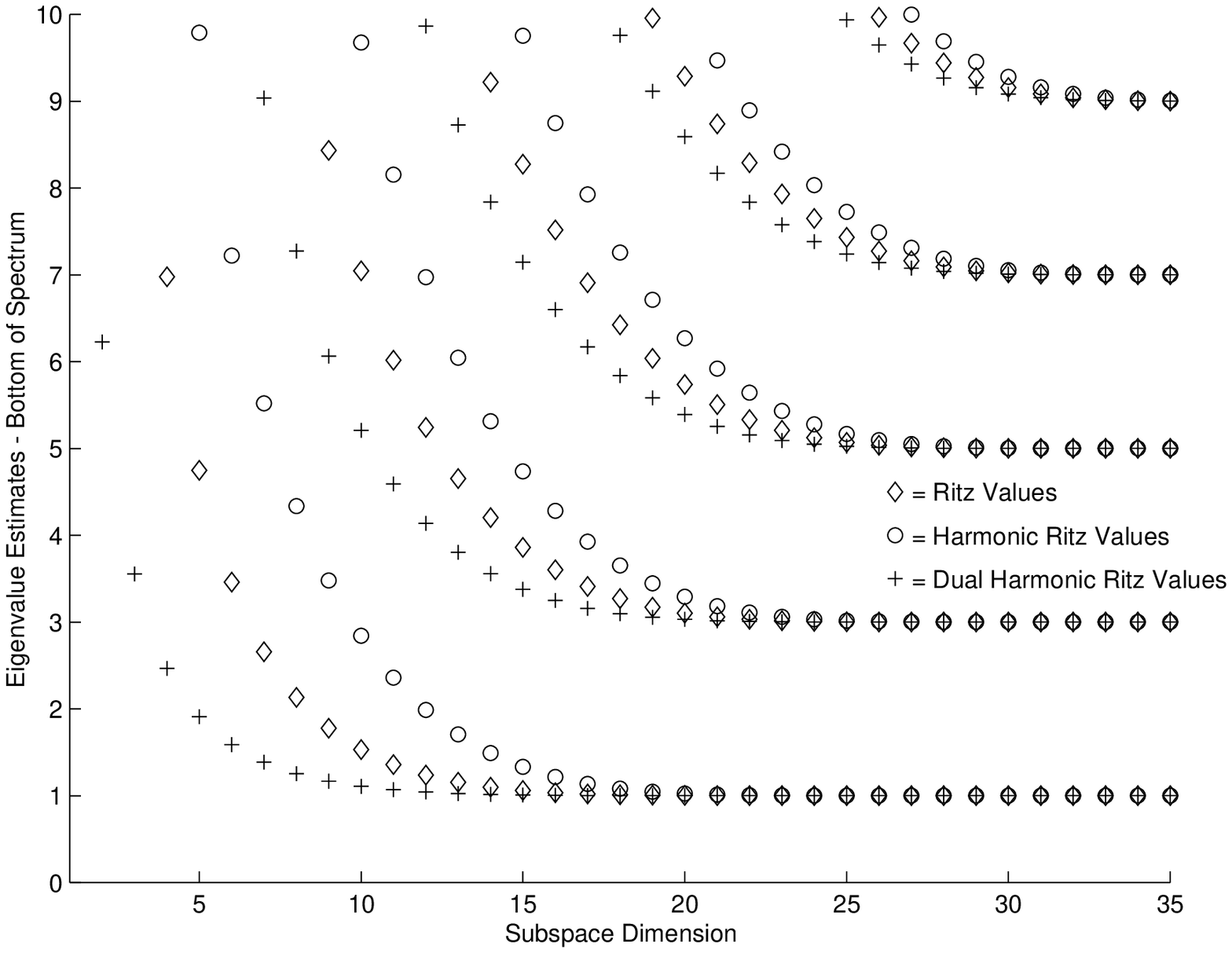,height=3in}}
 \caption{ Comparison of bounds on upper and  
		lower portions of the  spectrum.  }
 \label{fig:Ritz1}
 \end{figure}
 
 Both harmonic Ritz and dual harmonic 
Ritz values were known even 50 years ago and found to be useful in 
differential eigenvalue problems  --- Collatz 
\cite{Coll1}
referred to the harmonic Ritz problem (\ref{HarmRitzeig}) as Grammel's 
equations (citing Grammel's earlier work \cite{Grammel}) and viewed 
the Rayleigh quotients for the Ritz problem (\ref{Ritzeig}), the 
harmonic Ritz problem (\ref{HarmRitzeig}), and the dual harmonic Ritz 
problem (\ref{DualHarmRitzeig}), all as elements of an infinite 
monotone sequence of ``Schwarz quotients'' that could be generated 
iteratively.

As long as $\bK$ and $\bM$ are positive definite, 
all three of  Ritz, harmonic Ritz, and dual harmonic Ritz values provide ``inner'' bounds to the 
``outer'' eigenvalues of the pencil $\bK-\lambda \bM$ 
(that is, of the problem (\ref{eigprob})).  In comparing  the three 
types of approximations using the same subspace ${\mathcal P}$, 
harmonic Ritz values provide the best bounds of the three to the upper 
eigenvalues of (\ref{eigprob}); dual harmonic Ritz values provide the best bounds 
of the three to the lower eigenvalues.  As an example, Figure 1 shows 
bounds obtained for a sequence of nested Krylov subspaces taken for 
${\mathcal P}$, with $\mathbf{K}=diag([1:2:100])$, 
$\mathbf{M}=\mathbf{I}$, and a starting vector of all ones (the 
example of \cite{PaiParVor}).

\begin{theorem} \label{RitzBounds}
	Suppose  $\bK$ and $\bM$ are positive definite.  Then
		\begin{eqnarray*}
\lambda_{k}\ \leq\ \DRLambda_{k} \ \leq\ \Lambda_{k}\ \leq\ 
\tilde{\Lambda}_{k} & \mbox{for} & k=1,\ 2,\ \dots \\
     \DRLambda_{-\ell} \ \leq\ \Lambda_{-\ell}\ \leq\ 
     \tilde{\Lambda}_{-\ell}\ \leq\ \lambda_{-\ell} & \mbox{for} & \ell 
     =1,\ 2,\ \dots
	\end{eqnarray*}
\end{theorem}
\textsc{ Proof:} The min-max characterization 
yields
\begin{align*}
	\lambda_{k} = \min_{\dim {\mathcal S} = k}\max_{\bx\in {\mathcal S}}
	\frac{\bx^{t}\bK\bx}{\bx^{t}\bM\bx} & \leq
	\min_{\substack{
	\dim {\mathcal S} = k \\
	{\mathcal S}\subset {\mathcal P}
	}}\max_{\bx\in {\mathcal S}}
	\frac{\bx^{t}\bK\bx}{\bx^{t}\bM\bx} \\
	& = \min_{\dim {\mathcal R} = k}\max_{\by\in {\mathcal R}}
	\frac{\by^{t}\bP^{t}\bK\bP \by}{\by^{t}\bP^{t}\bM\bP \by}=\Lambda_{k},
\end{align*}
and likewise,
\begin{align*}
	\lambda_{k} = \min_{\dim {\mathcal S} = k}\max_{\bx\in {\mathcal S}}
	\frac{\bx^{t}\bK	\bM^{-1} \bK\bx}{\bx^{t}\bK\bx} & \leq
	\min_{\substack{
	\dim {\mathcal S} = k \\
	{\mathcal S}\subset {\mathcal P}
	}}\max_{\bx\in {\mathcal S}}
	\frac{\bx^{t}\bK	\bM^{-1} \bK\bx}{\bx^{t}\bK\bx} \\
	& = \min_{\dim {\mathcal R} = k}\max_{\by\in {\mathcal R}}
	\frac{\by^{t}\bP^{t}\bK	\bM^{-1} \bK\bP \by}{\by^{t}\bP^{t}\bK\bP \by}=\tilde{\Lambda}_{k}.
\end{align*}
A similar argument shows $\lambda_{k}\leq \DRLambda_{k}$.  By 
repeating the argument for the eigenvalue problem 
$-\bK\bx=(-\lambda)\bM\bx$, one finds $-\lambda_{\ell}(-\bK,\ \bM) 
\leq -\Lambda_{-\ell}$ (where $\lambda(\bA,\ \bB)$ is used to denote 
an eigenvalue of the pencil $\bA-\lambda \bB$).  Notice that 
$-\lambda_{\ell}(-\bK,\ \bM)=\lambda_{-\ell}(\bK,\ \bM)$.  Thus, 
$\Lambda_{-\ell}\leq \lambda_{-\ell}$ and $\tilde{\Lambda}_{-\ell}\leq 
\lambda_{-\ell}$.

For any $\bx \in \mathbb{R}^{n}$, the Cauchy-Schwarz 
inequality implies
\begin{align*}	
	(\bx^{t}\bK\bx)^{2}=& (\bx^{t}\bK\bM^{-1/2}\bM^{1/2}\bx)^{2}\leq 
	\bx^{t}\bK\bM^{-1}\bK\bx 	\ \bx^{t}\bM \bx \\
	\mbox{and} \quad (\bx^{t}\bM\bx)^{2}=& (\bx^{t}\bM\bK^{-1/2}\bK^{1/2}\bx)^{2}\leq 
	\bx^{t}\bM\bK^{-1}\bM\bx 	\ \bx^{t}\bK \bx 
\end{align*}

Thus, 
$$
\frac{\bx^{t}\bM\bx}{\bx^{t}\bM\bK^{-1}\bM\bx} \leq
\frac{\bx^{t}\bK\bx}{\bx^{t}\bM \bx}\leq 
\frac{\bx^{t}\bK\bM^{-1}\bK\bx}{\bx^{t}\bK\bx},
$$ which then implies
for each $k=1,\ 2,\ \dots,\ m$
\begin{align*}
0<\lambda_{k}\leq \DRLambda_{k} & = \min_{\substack{ \dim {\mathcal S} 
= k \\
		{\mathcal S}\subset {\mathcal P}
		}}\max_{\bx\in {\mathcal S}}
		\frac{\bx^{t}\bM\bx}{\bx^{t}\bM\bK^{-1}\bM\bx} \\
		& \leq  \min_{\substack{
		\dim {\mathcal S} = k \\
		{\mathcal S}\subset {\mathcal P}
		}}\max_{\bx\in {\mathcal S}}
		\frac{\bx^{t}\bK\bx}{\bx^{t}\bM\bx} = \Lambda_{k}\\
		& \leq \min_{\substack{
		\dim {\mathcal S} = k \\
		{\mathcal S}\subset {\mathcal P}
		}}\max_{\bx\in {\mathcal S}}
		\frac{\bx^{t}\bK \bM^{-1} \bK\bx}{\bx^{t}\bK\bx} =\tilde{\Lambda}_{k} \quad \blacksquare
\end{align*}

The situation is somewhat different if $\bK$ is \emph{indefinite}.
The Ritz estimates are still ``inner'' bounds, that is
$
\lambda_{k}\leq \Lambda_{k} $ and $ \Lambda_{-\ell} \leq\lambda_{-\ell}.
$
However, both harmonic Ritz and dual harmonic Ritz values  now provide
``outer''  bounds (\emph{lower} bounds) to negative eigenvalues of (\ref{eigprob})
and no simple relationship is known that would predict which of the 
three bounds is best (essentially owing to there being  no simple analog 
of the Cauchy-Schwarz inequality for indefinite inner products).

Despite  the differences in  behavior described above, Ritz, harmonic 
Ritz, and dual harmonic Ritz values each provide \emph{optimal} bounds -- obviously each with 
respect to a slightly different notion of optimality.
For the Ritz problem, the matrices $\bP^{t}\bK\bP$ and $\bP^{t}\bM\bP$ provide a ``sampling'' 
of the full matrices $\bK$ and $\bM$ on the subspace ${\mathcal P}$.  Whatever spectral 
information about the original eigenvalue problem (\ref{eigprob}) that we are able to deduce by 
examining the Rayleigh-Ritz problem (\ref{Ritzeig}) we must draw the 
same conclusions for 
{\em all}  matrix pencils that are ``aliased'' by the 
Rayleigh-Ritz sampling.  Define the following set of such $n \times n$ 
matrix pairs:
$$
{\mathcal C}({\mathcal P})=\left\{ (\bA,\ \bB) \left| 
\begin{array}{c}
\mbox{$\bA$ and $\bB$ are positive definite} \\
	\bP^{t}(\bA-\bK)\bP=0 \\
\bP^{t}(\bB-\bM)\bP=0  
\end{array}\right. \right\}
$$

\begin{theorem}
	For any choice of positive integers $\nu,\ \pi$ with $\nu+\pi=m$ and 
	any choice of matrix pairs $(\bA,\ \bB)\in {\mathcal C}({\mathcal P})$
	\begin{eqnarray*}
		\lambda_{k}(\bA,\ \bB) \leq \Lambda_{k}& \mbox{for}  & k=1,\ 2,\ 
		\dots,\ \nu  \\
	\Lambda_{-\ell}	 	\leq \lambda_{-\ell}(\bA,\ \bB) &  \mbox{for}  & 
\ell =1,\ 2,\ \dots,\ \pi.
	\end{eqnarray*}
Furthermore, for each index pair $\nu$, $\pi$, there exists a matrix pair 
$(\hat{\bA},\ \hat{\bB})\in {\mathcal C}({\mathcal P})$ such that 
	\begin{eqnarray*}
		\lambda_{k}(\hat{\bA},\ \hat{\bB}) = \Lambda_{k}& \mbox{for}  & k=1,\ 2,\ 
		\dots,\ \nu  \\
	\Lambda_{-\ell}	 	= \lambda_{-\ell}(\hat{\bA},\ \hat{\bB}) &  \mbox{for}  & 
\ell =1,\ 2,\ \dots,\ \pi.
	\end{eqnarray*}
  So, no better bounds 
are possible with only the information available to the Rayleigh-Ritz 
method as described by (\ref{Ritzeig}).
\end{theorem}
\textsc{ Proof:} The first assertion is a restatement of Theorem \ref{RitzBounds}
for the matrix pencil $\bA-\lambda \bB$.
 To show optimality, define the matrix of Ritz vectors:
$\bU=\left[\bu_{1} ,\
 \bu_{2} ,\ \dots, \bu_{\nu},\  \bu_{-\pi},\ \dots,\
 \bu_{-2},\ \bu_{-1}\right]$. 
 Notice that $\bU$ is an $\bM$-orthonormal basis for ${\mathcal P}$: 
 $\bU^{t}\bM\bU=\bI$.
 Define also 
 the diagonal matrix of Ritz values
 $$
 \bD=\left[
  \begin{array}{cccccc}
  	 \Lambda_{1} &  &  &  &  &   \\
  	 & \ddots &  &  &  &   \\
  	 &  & \Lambda_{\nu} &  &  &   \\
  	 &  &  & \Lambda_{-\pi} &  &   \\
  	 &  &  &  & \ddots &   \\
  	 &  &  &  &  & \Lambda_{-1}
  \end{array}\right]
 $$
 and fix $\hat{\Lambda}=\frac{1}{2}(\Lambda_{\nu}+\Lambda_{-\pi})$.
 Now, consider
$$
	 \hat{\bA} =\bM\bU\bD\bU^{t}\bM+ 
	 \hat{\Lambda}(\bM-\bM\bU\bU^{t}\bM)
	 \quad\mbox{and}\quad	
	 \hat{\bB} =\bM.
$$
One may verify that all required conditions are satisfied, in 
 particular
 $$(\hat{\bA}-\lambda\hat{\bB})\bU=\bM\bU(\bD-\lambda\bI)$$
 and for any
 $\bv\in \mathbb{R}^{n}$ with $\bv^{t}\bM\bU=0$,
 $$(\hat{\bA}-\lambda\hat{\bB})\bv=\bM\bv(\hat{\Lambda}-\lambda).
 \quad \blacksquare
$$
A similar construction can be used to show the (analogously defined) 
optimality of harmonic Ritz values and dual harmonic Ritz values.

As we will see in following sections, Ritz values, harmonic Ritz values, 
and dual harmonic Ritz values are limiting cases of parameterized 
families of bounds arising from ``left-definite'' and ``right-definite'' Lehmann 
intervals.

\setlength{\unitlength}{1.3mm}
\begin{picture}(85,40)(-7,-6)
\thinlines \put(0,0){\line(4,3){8.5}}
           \put(19.6,14.7){\line(4,3){16.4}} 
           \put(36,27){\line(4,-3){16.4}} 
           \put(63.5,6.375){\line(4,-3){8.5}} 
\put(-6,-4){\makebox(0,0){\begin{tabular}{|c|} \hline
	 Dual Harmonic \\
	 Ritz Values \\
	\hline
\end{tabular} }}
\put(19.8,10){\oval(25.5,9)}
\put(20,10){\makebox(0,0){\begin{tabular}{c}
Left-definite \\
Lehmann Bounds\\
\end{tabular} }}
\put(36,29){\makebox(0,0){\begin{tabular}{|c|}
	\hline
	 Ritz Values \\
	\hline
\end{tabular}}}
\put(25,24){\makebox(0,0){$\rho\rightarrow 0$}}
\put(48,24){\makebox(0,0){$\rho\rightarrow \pm \infty$}}
\put(76,3){\makebox(0,0){$\rho\rightarrow 0$}}
\put(-3,3){\makebox(0,0){$\rho\rightarrow \pm \infty$}}
\put(79,-4){\makebox(0,0){\begin{tabular}{|c|}
	\hline
	 Harmonic \\
	 Ritz Values \\
	\hline
\end{tabular}}}
\put(51.8,10){\oval(25.5,9)}
\put(52,10){\makebox(0,0){\begin{tabular}{c}
Right-definite \\
Lehmann bounds\\
\end{tabular} }}
\end{picture}

\section{Lehmann's Optimal Intervals}
Each of the Ritz-related methods discussed above will have certain advantages in estimating 
the extreme eigenvalues of (\ref{eigprob}).  None are particularly 
effective in estimating interior eigenvalues, however.  Usual 
strategies for obtaining accurate estimates to the eigenvalues of 
(\ref{eigprob}) lying close to a given value $\rho$ involve a spectral 
mapping that turns the spectrum ``inside out'' around $\rho$ --- 
mapping interior eigenvalues in the neighborhood of $\rho$ to extreme 
eigenvalues that are more accessible.  ``Shift and invert'' strategies 
typically use the spectral mapping $\lambda \mapsto \frac{1}{\lambda - 
\rho}$.  A variant used especially for buckling problems (where $\bM$ 
may be singular) utilizes instead the spectral mapping $\lambda \mapsto 
\frac{\lambda}{\lambda - \rho}$.  As we shall see, both of these spectral mappings play 
a fundamental role in the optimal bounds discovered by Lehmann 
(\cite{Leh49}, \cite{Leh50}, \cite{Leh63}).  The derivation used here 
is in the spirit of that given by Maehly in \cite{Maehly} and the 
associated methods are sometimes called Lehmann-Maehly methods.

Fix a scalar
$\rho$ that is not an eigenvalue of (\ref{eigprob}) and define the index $r$ 
to satisfy
\begin{equation}
	\lambda_{r-1} < \rho < \lambda_{r}.
	\label{seppar}
\end{equation}

The \emph{right-definite Lehmann method} follows first from considering the spectral mapping 
$\lambda \mapsto \frac{1}{\lambda - \rho}$ and an associated eigenvalue problem 
equivalent to (\ref{eigprob}):
\begin{equation}
	\bM(\bK-\rho\bM)^{-1}\bM\bx = \frac{1}{\lambda - \rho} \bM\bx,
	\label{Rdefeig}
\end{equation}
which has eigenvalues distributed as
\begin{equation*}
	\frac{1}{\lambda_{r-1} - \rho} \leq \frac{1}{\lambda_{r-2} - \rho} 
	\leq \dots \leq 0 \leq  \dots \leq
	\frac{1}{\lambda_{r+1} - \rho} \leq \frac{1}{\lambda_{r} - \rho}
\end{equation*}
Notice that eigenvalues of (\ref{eigprob}) flanking $\rho$ are mapped 
to extremal eigenvalues of (\ref{Rdefeig}).  Now use an 
$m$-dimensional subspace ${\mathcal S}=Ran(\bS)$ to generate Rayleigh-Ritz estimates for 
the eigenvalues of (\ref{Rdefeig}):
\begin{equation}
	[\bS^{t}\bM(\bK-\rho\bM)^{-1}\bM\bS] \by = R \ [\bS^{t}\bM\bS] \by,
	\label{RdefRitz}
\end{equation}
where $\bS\in \mathbb{R}^{n\times m}$.
Suppose (\ref{RdefRitz}) has $\nu$ negative eigenvalues
$R_{1}\leq \dots \leq R_{\nu}<0 $ and $\pi=m-\nu$ positive 
eigenvalues $0<R_{-\pi}\leq \dots \leq R_{-1} $. 
Regardless of the subspace ${\mathcal S}$ that is chosen, the 
min-max principle (or Theorem \ref{RitzBounds}) guarantees that
for each $k=1,\ 2,\ \dots, \nu$ and $\ell=1,\ 2,\ \dots, \pi$
\begin{eqnarray*}
\frac{1}{\lambda_{r-k} - \rho} \leq  R_{k}	 &\mbox{ and }  &
 R_{-\ell}\leq \frac{1}{\lambda_{r+\ell-1} - \rho}.
\end{eqnarray*}
Rearrange and introduce
\begin{eqnarray} \label{Rdefbounds}
 \Lambda^{(R)}_{-k}\stackrel{def}{=}\rho+\frac{1}{  R_{k} } \leq \lambda_{r-k}	 &\mbox{ and }  &
\lambda_{r+\ell-1} \leq \rho+\frac{1}{  R_{-\ell} 
}\stackrel{def}{=}\Lambda^{(R)}_{\ell}
\end{eqnarray}
for $k=1,\ 2,\ \dots, \nu$ and $\ell=1,\ 2,\ \dots, \pi$. Notice that 
labeling of $\Lambda^{(R)}$ is arranged relative to $\rho$:
$$
 \dots\Lambda^{(R)}_{-3}\leq \Lambda^{(R)}_{-2} \leq \Lambda^{(R)}_{-1}  
 < \rho < \Lambda^{(R)}_{1} 
\leq \Lambda^{(R)}_{2} \leq \Lambda^{(R)}_{3}\dots
$$ 
An equivalent statement combining (\ref{seppar})  and (\ref{Rdefbounds}) is \\
\begin{center}
	\begin{minipage}[c]{4.5in}
		Each of the intervals $[ \Lambda^{(R)}_{-k},\, \rho)$ and  $(\rho,\,  
		\Lambda^{(R)}_{\ell}]$ contain respectively $k$ and $\ell$ eigenvalues 
		of (\ref{eigprob}) for $k=1,\ 2,\ \dots, \nu$ and $\ell=1,\ 2,\ \dots, \pi$.
	\end{minipage}
\end{center}
\vskip .4cm

To avoid the need in (\ref{RdefRitz}) for solving linear systems having the indefinite 
coefficient matrix  $(\bK-\rho\bM)$, change 
variables in (\ref{RdefRitz}) as $\bP=(\bK-\rho\bM)^{-1}\bM\bS$ ---
which then \emph{implicitly} determines ${\mathcal S}$ via a choice of 
${\mathcal P}$.
(\ref{RdefRitz}) can then be rewritten as
\begin{equation}
	[\bP^{t}(\bK-\rho\bM)\bP] \by = R\ [\bP^{t}(\bK-\rho\bM)\bM^{-1}(\bK-\rho\bM)\bP]\by 
	\label{eq:RdefLehm}
\end{equation}
When $\dim {\mathcal P}=1$, (\ref{Rdefbounds}) becomes \emph{Temple's 
inequality}
$$
\rho+\frac{\bp^{t}(\bK-\rho\bM)\bM^{-1}(\bK-\rho\bM)\bp}{ \bp^{t}(\bK-\rho\bM)\bp  } 
=\frac{\bp^{t}(\bK\bM^{-1}\bK-\rho\bK)\bp}{ \bp^{t}(\bK-\rho\bM)\bp  } \leq \lambda_{r-1}
$$

Some additional notation will reduce the impending clutter of 
symbols.  Introduce matrices of \emph{Schwarz constants}:
$$
\bH_{0}=[\bP^{t}\bK\bM^{-1}\bK\bP],\
\bH_{1}=[\bP^{t}\bK\bP],\ \mbox{ and }\ 
\bH_{2}=[\bP^{t}\bM\bP].
$$
Then expanding out the various terms, (\ref{eq:RdefLehm}) becomes
\begin{equation}
	[\bH_{1}-\rho \bH_{2}]\by = R \, [\bH_{0}-2 \rho \bH_{1} + 
	\rho^{2}\bH_{2}]\by 
	\label{RdefPform}
\end{equation}
which may be rearranged to obtain
\begin{equation}
	[\bH_{0}- \rho \bH_{1}]\by = \Lambda^{(R)} \, [\bH_{1}-\rho \bH_{2}] \by 
	\label{RdefJform}
\end{equation}

Notice that (\ref{RdefJform}) could be written in terms of the $\bM$-inner product as
\begin{equation}
\bP^{\mathsf{m}}[(\bM^{-1}\bK)^{2}- \rho (\bM^{-1}\bK)]\bP\by = \Lambda^{(R)} \, 
\bP^{\mathsf{m}}[(\bM^{-1}\bK)-\rho \bI]\bP \by
\label{RdefMgenform}
\end{equation}
or in terms of the $\bK$-inner product as
\begin{equation}
\bP^{\mathsf{k}}[(\bM^{-1}\bK)- \rho \bI]\bP\by = \Lambda^{(R)} \, 
\bP^{\mathsf{k}}[\bI-\rho (\bM^{-1}\bK)^{-1}]\bP \by
\label{RdefKgenform}
\end{equation}
\vskip .5cm

The \emph{left-definite Lehmann method} can be obtained by 
considering the spectral mapping $\lambda \mapsto 
\frac{\lambda}{\lambda - \rho}$ and an associated eigenvalue problem 
--- also equivalent to (\ref{eigprob}):
\begin{equation}
	\bK(\bK-\rho\bM)^{-1}\bK\bx = \frac{\lambda}{\lambda - \rho} \bK\bx 
	\label{Ldefeig}
\end{equation}
which has eigenvalues distributed as
\begin{equation}
	\frac{\lambda_{r-1}}{\lambda_{r-1} - \rho} \leq \frac{\lambda_{r-2}}{\lambda_{r-2} - \rho} 
	\leq \dots < 0 \mbox{ and } 1 < \dots \leq   
	\frac{\lambda_{r+1}}{\lambda_{r+1} - \rho} \leq \frac{\lambda_{r}}{\lambda_{r} - \rho}
	\label{ldefspec}
\end{equation}
(as long as both $\bK$ and $\bM$ are positive definite, no eigenvalue 
gets mapped into the interval $[0,1]$).
Again the eigenvalues of (\ref{eigprob}) flanking $\rho$ are mapped 
to extremal eigenvalues of (\ref{Ldefeig}).  Using an 
$m$-dimensional subspace ${\mathcal T}=Ran(\bT)$,
one may generate Rayleigh-Ritz estimates for 
the eigenvalues of (\ref{Ldefeig}):
\begin{equation}
	[\bT^{t}\bK(\bK-\rho\bM)^{-1}\bK\bT] \by = L [\bT^{t}\bK\bT] \by,
	\label{LdefRitz}
\end{equation}
where $\bT\in \mathbb{R}^{n\times m}$.

If (\ref{LdefRitz}) has $\nu$ negative eigenvalues
$L_{1}\leq L_{2}\leq \dots \leq L_{\nu}<0 $ and $\pi=m-\nu$ positive 
eigenvalues $1<L_{-\pi}\leq \dots \leq L_{-2} \leq L_{-1} $, then 
regardless of the subspace ${\mathcal T}$ that is chosen, the 
min-max principle (or again, Theorem \ref{RitzBounds}) guarantees that
\begin{eqnarray}\label{LLbnds}
\frac{\lambda_{r-k}}{\lambda_{r-k} - \rho} \leq  L_{k}	 &\mbox{ and }  &
 L_{-\ell}\leq \frac{\lambda_{r+\ell-1}}{\lambda_{r+\ell-1} - \rho}  
\end{eqnarray}
or equivalently
\begin{eqnarray} \label{Ldefbounds}
 \quad\quad \Lambda^{(L)}_{-k}\stackrel{def}{=} \rho-\frac{\rho}{1- L_{k} 
 } \leq \lambda_{r-k} &\mbox{ and } & \lambda_{r+\ell-1} \leq 
 \rho-\frac{\rho}{1 - L_{-\ell} } \stackrel{def}{=}\Lambda^{(L)}_{\ell}
\end{eqnarray}
for $k=1,\ 2,\ \dots, \nu$ and $\ell=1,\ 2,\ \dots, \pi$.
Just as for $\Lambda^{(R)}$,  the labeling of  $\Lambda^{(L)}$ is done relative to $\rho$:
$$
 \dots\Lambda^{(L)}_{-3}\leq \Lambda^{(L)}_{-2} \leq \Lambda^{(L)}_{-1}  
 < \rho < \Lambda^{(L)}_{1} 
\leq \Lambda^{(L)}_{2} \leq \Lambda^{(L)}_{3}\dots
$$ 
An equivalent statement combining (\ref{seppar}) and (\ref{Ldefbounds}) 
is  \\
\begin{center}
	\begin{minipage}[c]{4.5in}
		Each of	the intervals $[ \Lambda^{(L)}_{-k},\, \rho)$ and  $(\rho,\,  
		\Lambda^{(L)}_{\ell}]$ contain respectively $k$ and $\ell$ eigenvalues 
		of (\ref{eigprob}) for  $k=1,\ 2,\ \dots, \nu$ and $\ell=1,\ 2,\ \dots, \pi$.
	\end{minipage}
\end{center}
\vskip .4cm

As before, in order to avoid solving systems with the coefficient 
matrix $(\bK-\rho\bM)$, change 
variables in (\ref{LdefRitz}) as $\bP=(\bK-\rho\bM)^{-1}\bK\bT$
which then \emph{implicitly} determines ${\mathcal T}$ via a choice of 
${\mathcal P}$.
(\ref{LdefRitz})  can then be rewritten as
\begin{equation}
	[\bP^{t}(\bK-\rho\bM)\bP] \by = L\ [\bP^{t}(\bK-\rho\bM)\bK^{-1}(\bK-\rho\bM)\bP]\by.
	\label{eq:LdefLehm}
\end{equation}

Introduce
$$
\bH_{3}=[\bP^{t}\bM\bK^{-1}\bM\bP].
$$
Then (\ref{eq:LdefLehm}) becomes
\begin{equation}
	[\bH_{1}-\rho \bH_{2}]\by = L [\bH_{1}-2 \rho \bH_{2} + 
	\rho^{2}\bH_{3}]\by.
	\label{LdefPform}
\end{equation}
which may be rearranged to get
 \begin{equation}
	[\bH_{1}-\rho \bH_{2}]\by = \Lambda^{(L)}\, [\bH_{2} -	\rho \bH_{3}]\by.
	\label{LdefJform}
\end{equation}
\vskip .5cm

Observe that both (\ref{RdefPform}) and 
(\ref{LdefPform}) are Hermitian definite pencils with the same 
left-hand side.  By the Sylvester Law of Inertia, they each have the 
same number of negative (and hence positive) eigenvalues.
If a shift of $\rho=0$ is chosen in (\ref{RdefJform}), the harmonic 
Ritz problem (\ref{HarmRitzeig}) is obtained and 
$\tilde{\Lambda}_{\ell}=\Lambda^{(R)}_{\ell}\left|_{\rho=0}\right.$.  
As $\rho\rightarrow \pm\infty$, (\ref{RdefJform}) reduces to the Ritz 
problem (\ref{Ritzeig}).  Similarly, if a shift of $\rho=0$ is chosen in 
(\ref{LdefJform}), the  Ritz problem (\ref{Ritzeig}) 
is obtained and $\Lambda_{\ell}=\Lambda^{(L)}_{\ell}\left|_{\rho=0}\right.$.  
As $\rho\rightarrow \pm\infty$, (\ref{LdefJform}) reduces to the dual 
harmonic Ritz problem (\ref{DualHarmRitzeig}).

The left- and right-definite Lehmann bounds, $\Lambda^{(L)}$ and 
$\Lambda^{(R)}$, that are below the parameter $\rho$ are monotone increasing 
with respect to  $\rho$.  This is easy to show for $\rho$ satisfying 
(\ref{seppar}), however as $\rho$ is increased further, $r$ changes 
and the labeling of $\Lambda^{(L)}$ and $\Lambda^{(R)}$ shifts.  This 
more complicated circumstance is dicussed in \cite{Zimm} where a proof 
of monotonicity in the general case may be found.  
 
Notice that (\ref{LdefJform}) could be obtained formally from the 
right-definite method expressed in (\ref{RdefKgenform}) by direct 
substitution of the $\bM$-inner product for the $\bK$-inner product.
\begin{equation}
\bP^{\mathsf{m}}[(\bM^{-1}\bK)- \rho \bI]\bP\by = \Lambda^{(L)} \, 
\bP^{\mathsf{m}}[\bI-\rho (\bM^{-1}\bK)^{-1}]\bP \by
\label{LdefMgenform}
\end{equation}
Such a substitution also converts the harmonic Ritz problem into a Ritz 
problem and the Ritz problem, then into a dual harmonic Ritz problem.  
This provides some impetus to call the ``left-definite Lehmann'' 
method the ``harmonic Lehmann'' method, but Lehmann himself referred 
to this method as ``left-definite'' and besides the correspondences 
are a bit backward since (right-definite) Lehmann is to Ritz as ``dual 
harmonic Ritz'' is to `` harmonic Lehmann.''

\section{An Alternative Formulation}
Kahan developed a formulation of Lehmann's right-definite method that is 
particularly well-suited to many computational settings for matrix 
eigenvalue problems (cf.  \cite{Parlett}, Chap. 10).  We review the 
development here and extend it to Lehmann's left-definite method.  For a 
given $m$-dimensional subspace ${\mathcal P}$, suppose the columns of 
$\bQ_{1}$ provide an $\bM$-orthonormal basis for ${\mathcal P}$: 
$span(\bQ_{1})={\mathcal P}$ and 
$\bQ_{1}^{\mathsf{m}}\bQ_{1}=\bQ_{1}^{t}\bM\bQ_{1}=\bI$.  Define $\bH$ 
from the ``residual orthogonality'' condition
$$
(\bM^{-1}\bK\bQ_{1}-\bQ_{1}\bH)^{t}\bM\bQ_{1}=0
$$
so that $\bH=\bQ_{1}^{t}\bK\bQ_{1}$ and observe (say, from the 
Gram-Schmidt process) that there is an upper 
triangular matrix $\bC$ and a matrix $\bQ_{2}$ with $\bM$-orthonormal 
columns so that 
$$
\bQ_{2}\bC=\bM^{-1}\bK\bQ_{1}-\bQ_{1}\bH.
$$
Pick $\bQ_{3}$ to fill out an $\bM$-orthonormal basis for ${\mathbb 
R}^{n}$ in conjunction with $\bQ_{1}$ and $\bQ_{2}$.  Then with 
$\bQ=\left[\bQ_1\ \bQ_2\ \bQ_3\right]$, we have $\bQ^{t}\bM\bQ=\bI$ 
and
$$
\bM^{-1}\bK\bQ=\bQ\left[
\begin{array}{ccc}
	 \bH & \bC^{t} & 0  \\
	\bC & \bV_{11} &  \bV_{21}^{t} \\
	0 & \bV_{21} & \bV_{22}
\end{array}\right]\quad
\begin{array}{l}
	\mbox{where}  \\
	\ \bH \mbox{ is } m\times m \\
		\ \bV_{11} \mbox{ is } k \times k.
\end{array}
$$
While this shows how $\bH$ and $\bC$ might be constructed (essentially 
one step of a block Lanczos process), there may be other situations of 
interest when $\bH$ and $\bC$ are known \emph{a priori}.  In any case, 
we assume that the bottom right block $2\times 2$ submatrix, $\bV$, is 
either unknown or at least unpleasant to deal with.  With additional 
unitary massage, $rank(\bC)=k$ could be assumed (possibly resulting in 
a smaller $\bV_{11}$), though it isn't necessary in what follows.  The 
situation $rank(\bC) = k \ll m \ll n$ is common.  What follows is a  
\emph{deus ex machina} development of Kahan's formulation of Lehmann 
bounds  that offers brevity but little of the insight 
and revelation that one may find in the excellent discussion of 
(\cite{Parlett}, Chapter 10).

Apply the right-definite Lehmann bounds from (\ref{eq:RdefLehm}) using 
$\bP=\bQ_1$.  Then, $(\bK-\rho\bM)\bP=\bQ_1(\bH-\rho\bI)+\bQ_2\bC$ and 
the right-definite Lehmann problem (\ref{RdefPform}) appears as
\begin{equation}
	(\bH-\rho\bI)\by=R\ 
	\left[(\bH-\rho\bI)^{2}+\bC^{t}\bC\right]\by
	\label{eq:kahan1}
\end{equation}
The associated  right-definite bound is $\Lambda^{(R)}=\rho+1/R$ and we
may manipulate (\ref{eq:kahan1}) to 
get an equivalent condition on  
$\Lambda^{(R)}$:
\begin{equation} 
	0=\left[(\bH-\rho\bI)(\bH-\Lambda^{(R)}\bI)+\bC^{t}\bC\right]\by
	\label{eq:kahan2}
\end{equation}
One may recognize that the coefficient matrix of (\ref{eq:kahan2}) is 
a Schur complement of the $(m+k)\times (m+k)$ matrix 
$$\bY(\Lambda^{(R)})=\left[
\begin{array}{cc}
	-(\bH-\rho\bI)(\bH-\Lambda^{(R)}\bI) & \bC^{t}  \\
	\bC & \bI
\end{array}\right].
$$
Hence, (\ref{eq:kahan2}) has a non-trivial solution if and only if 
$\bY(\Lambda^{(R)})$ is singular.  Suppose that neither $\rho$ nor $\Lambda^{(R)}$ are eigenvalues of
$\bH$ for the time being and define
\begin{align*}
	\bL_{1} & =\left[
	\begin{array}{cc}
		\bI & 0  \\
		\bC(\bH-\rho\bI)^{-1}(\bH-\Lambda^{(R)}\bI)^{-1} & \bI
	\end{array}\right],  \\	
	\bL_{2} & =\left[
	\begin{array}{cc}
		\bI & 0  \\
		\bC(\bH-\Lambda^{(R)}\bI)^{-1} & \bI
	\end{array}\right], \\
	\mbox{  and  }\quad
	\bD(\Lambda^{(R)}) & =\left[
	\begin{array}{cc}
		-(\bH-\rho\bI)^{-1} & 0  \\
		0 & (\rho-\Lambda^{(R)})\bI
	\end{array}\right].
\end{align*}

Then
$$
\bL_{2}\bD(\Lambda^{(R)})\bL_{1}\bY(\Lambda^{(R)})\bL_{1}^{t}\bL_{2}^{t}
=\left[
\begin{array}{cc}
	\bH-\Lambda^{(R)}\bI & \bC^{t}  \\
	\bC & \rho\bI+\bC(\bH-\rho\bI)^{-1}\bC^{t}-\Lambda^{(R)}\bI
\end{array}\right].
$$
Thus  $\Lambda^{(R)}$ is an eigenvalue of the $(m+k) \times (m+k)$ matrix
\begin{equation}
	\left[
	\begin{array}{cc}
		\bH & \bC^{t}  \\
		\bC & \rho\bI+\bC(\bH-\rho\bI)^{-1}\bC^{t}
	\end{array}\right].
\label{kahfin}
\end{equation}
if and only if either $\bD(\Lambda^{(R)})$ is singular or 
$\bY(\Lambda^{(R)})$ is singular, which is to say,
if and only if either $\Lambda^{(R)}$ is a right-definite 
Lehmann bound satisfying (\ref{eq:kahan2}) or $\Lambda^{(R)}=\rho$ 
(which will occur with multiplicity $k$).  A limiting argument can be 
mustered to handle the exceptional cases where either $\rho$ or 
$\Lambda^{(R)}$ are eigenvalues of $\bH$.  In situations where either the 
smaller eigenvalues of (\ref{eigprob}) are of interest or  $\|\bC\|$ 
is much smaller than $\|\bH\|$, finding the eigenvalues of 
(\ref{kahfin}) is likely to yield substantially more accurate results 
for $\Lambda^{(R)}$ then a direct attack on (\ref{eq:kahan1}). A similar 
formulation for left-definite Lehmann 
problems will be described below.

Consider the application of the left-definite problem 
(\ref{LdefPform}) with $\bP=\bQ_1$.  Note that 
$\bK\bQ_1=\bQ_1\bH+\bQ_2\bC$ implies that 
$$
\bK^{-1}\bQ_1=\bQ_1\bH^{-1}-\bK^{-1}\bQ_2\bC\bH^{-1}
$$ so then 
\begin{equation}
	\bQ_1^{t}\bK^{-1}\bQ_1=\bH^{-1}+\bH^{-1}\bC^{t}\bW\bC\bH^{-1}
\label{res}
\end{equation}
 where $\bW=\bQ_2^{t}\bK^{-1}\bQ_2$ has been introduced.
 (\ref{LdefPform}) becomes
\begin{equation}
	(\bH-\rho\bI)\by=L\ 
	\left[(\bH-\rho\bI)-
\rho(\bI-\rho(\bH^{-1}+\bH^{-1}\bC^{t}\bW\bC\bH^{-1}))\right]\by.
	\label{eq:kahan3}
\end{equation}
The associated  left-definite bound is $\Lambda^{(L)}=-\rho L/(1-L)$ and we
may manipulate (\ref{eq:kahan3}) to 
get an equivalent condition on  
$\Lambda^{(L)}$:
\begin{equation} 
	0=\left[(\bH-\rho\bI)(\bH-\Lambda^{(L)}\bI)\bH+\rho \Lambda^{(L)} 
	\bC^{t}\bW\bC\right]\by.
	\label{eq:kahan4}
\end{equation}

(\ref{eq:kahan4}) has a non-trivial solution if and only if the 
$(m+k)\times (m+k)$ matrix
$$\hat{\bY}(\Lambda^{(L)})=\left[
\begin{array}{cc}
	-(\bH-\rho\bI)(\bH-\Lambda^{(L)}\bI)\bH & \Lambda^{(L)}\bC^{t}  \\
	\rho\bC & \bW^{-1}
\end{array}\right]
$$
is singular.  Suppose that neither $\rho$ nor $\Lambda^{(L)}$ are eigenvalues of
$\bH$ and define 
\begin{align*}
	\bF & =(\bH-\rho\bI)^{-1}(\bH-\Lambda^{(L)}\bI)^{-1}\bH^{-1}\\
	\hat{\bL}_{1}& =\left[
	\begin{array}{cc}
		\bI & 0  \\
\rho\bC\bF & \bI \end{array}\right], \quad 
\hat{\bU}_{1} =\left[
	\begin{array}{cc}
\bI & \Lambda^{(L)}\bF\bC^{t} \\
0 & \bI \end{array}\right],\\
	\hat{\bL}_{2}& =\left[
	\begin{array}{cc}
		\bI & 0  \\
		\bC(\bH-\Lambda^{(L)}\bI)^{-1} & \bI
	\end{array}\right], \\
	\mbox{  and  }\quad
	\hat{\bD}(\Lambda^{(L)}) & =\left[
	\begin{array}{cc}
		-(\bH-\rho\bI)^{-1}\bH^{-1} & 0  \\
		0 & (\rho-\Lambda^{(L)})/\rho\bI
	\end{array}\right].
\end{align*}

Then
	\begin{equation}																		 
		\hat{\bL}_{2}\hat{\bD}(\Lambda^{(L)})\hat{\bL}_{1}\hat{\bY}(\Lambda^{(L)})			 
	\hat{\bU}_{1}\hat{\bL}_{2}^{t} =\left[													 
		\begin{array}{cc}																	 
			\bH-\Lambda^{(L)}\bI & \bC^{t}	\\												 
			\bC	& \frac{\rho-\Lambda^{(L)}}{\rho}\bN_{1}+\frac{\Lambda^{(L)}}{\rho}\bN_{2}	 
		\end{array}\right].																	 
		\label{eq:kahan4_5}																	 
	\end{equation}		
where 
$\bN_{1}=\bW^{-1}+\bC\bH^{-1}\bC^{t}$ and $\bN_{2}=\bC(\bH-\rho\bI)^{-1}\bC^{t}$.
Thus  $\Lambda^{(L)}$ is an eigenvalue of an auxiliary $(m+k) \times (m+k)$ matrix 
pencil --- not unlike the right-definite case. 
This matrix pencil will be \emph{definite} when $\bN_{1}-\bN_{2}$ is 
positive-definite, which in turn can be guaranteed when the $(r-1)st$ 
Ritz value is a sufficiently accurate approximation to 
$\lambda_{r-1}$:
\begin{theorem} \label{leftdefKahan}
Suppose $\rho$ is not an eigenvalue of (\ref{eigprob}).  Each 
interval $[\Lambda^{(L)}_{-i},\ \rho)$ and $(\rho,\ 
\Lambda^{(L)}_{j}]$ contains respectively at least  $i$ and $j$ 
eigenvalues of (\ref{eigprob}), where
	$$
 0<\Lambda^{(L)}_{-\nu} \leq \dots \leq \Lambda^{(L)}_{-2} \leq 
 \Lambda^{(L)}_{-1}< \rho < \Lambda^{(L)}_1 \leq \Lambda^{(L)}_2 \leq 
 \dots
	$$
are the positive eigenvalues of the $(m+k) \times (m+k)$ matrix pencil
\begin{align}
	\left[
	\begin{array}{cc}
		\bH & \bC^{t}  \\
		\bC & \bN_{1}
	\end{array}\right]& - \Lambda^{(L)}\left[
	\begin{array}{cc}
		\bI & 0  \\
0 & \bM_{1} \end{array}\right], \label{eq:kahan5} \\
 \mbox{where }\quad \bM_{1}& =\frac{1}{\rho}(\bN_{1}-\bN_{2}) 
 \nonumber \\
 \bN_{1}& =\bW^{-1}+\bC\bH^{-1}\bC^{t},\quad \mbox{ and } \nonumber \\
	\bN_{2}& =\bC(\bH-\rho\bI)^{-1}\bC^{t}. \nonumber
\end{align}
$\rho$ is an eigenvalue of (\ref{eq:kahan5}) with multiplicity $k$.
If the Ritz value $\Lambda_{r-1}<\rho$, then $\bM_{1}$ is positive 
definite and (\ref{eq:kahan5}) is a Hermitian definite pencil.
\end{theorem}
\textsc{ Proof:} The first assertion follows immediately from 
(\ref{eq:kahan4_5}), since then $\Lambda^{(L)}$ is an eigenvalue of 
(\ref{eq:kahan5}) if and only if either $\hat{\bD}(\Lambda^{(L)})$ is 
singular or $\hat{\bY}(\Lambda^{(L)})$ is singular.  As before a 
limiting argument handles the exceptional cases where either $\rho$ or 
$\Lambda^{(L)}$ are eigenvalues of $\bH$.

For the second statement, note that $\Lambda_{r-1}<\rho$ implies from 
the way that $r$ was chosen in (\ref{seppar}) that
$\bH-\rho\bI$ has precisely $r-1$ negative eigenvalues.  Note then 
that $\bN_{1}-\bN_{2}$ is positive-definite  if and 
only if the matrix
\begin{equation}
	\left[
	\begin{array}{cc}
	  \frac{1}{\rho}(\bH-\rho\bI)\bH	 & 0  \\
		0 & \bN_{1}-\bN_{2}
	\end{array}\right]
	\label{eq:kahprf1}
\end{equation}
has precisely $r-1$ negative eigenvalues.  Define 
\begin{align*}
	\tilde{\bL}_{1}= \left[
	\begin{array}{cc}
		\bI & 0  \\
\rho\bC\bH^{-1}(\bH-\rho\bI)^{-1} & \bI \end{array}\right],&\\
	\tilde{\bL}_{2}=\left[
	\begin{array}{cc}
		\bI & - \bC^{t}\bW\\
		0& \bI
	\end{array}\right],  \quad \mbox{ and }\ \tilde{\bD}&=
	\left[
	\begin{array}{cc}
\rho\bH^{-1} & 0\\
		0 & \bI
	\end{array}\right].
\end{align*}
and calculate with $\tilde{\bF}=\tilde{\bD}\tilde{\bL}_{2}\tilde{\bL}_{1}$
\begin{equation}
	\tilde{\bF}\left[
		\begin{array}{cc}
\frac{1}{\rho}(\bH-\rho\bI)\bH & 0 \\
			0 & \bN_{1}-\bN_{2}
\end{array}\right]\tilde{\bF}^{t}= \left[
		\begin{array}{cc}
\rho(\bI-\rho \bP^{t}\bK^{-1}\bP) & 0 \\
		0	 & \bW^{-1}
		\end{array}\right].
	\label{eq:kahprf2}
\end{equation}
Suppose (\ref{eq:kahprf1}) had more than $r-1$ negative eigenvalues.
Then (\ref{eq:kahprf2}) has more than $r-1$ negative eigenvalues and 
therefore  $\bI-\rho \bP^{t}\bK^{-1}\bP $ has more than $r-1$ negative 
eigenvalues.  Equivalently, this means that $\bP^{t}\bK^{-1}\bP$ has $r$ 
or more eigenvalues \emph{above} $1/\rho$.  Since the eigenvalues of 
$\bP^{t}\bK^{-1}\bP$ provide inner bounds to the outer eigenvalues of 
$\bK^{-1}$, this implies in turn that $\bK^{-1}$ must have $r$ 
or more eigenvalues \emph{above} $1/\rho$.  But this contradicts the 
choice of $\rho$ made in (\ref{seppar}). $\quad \blacksquare$

The calculation of $\bW=\bQ_{2}^{t}\bK^{-1}\bQ_{2}$ involves the 
solution of $k$ linear systems each of the form $\bK\bx=\bb$.  If these 
systems are solved inexactly (one rarely has other options), 
reasonable concerns arise about the integrity of the resulting 
bounds.  Rigorous inclusion intervals can be maintained if the 
approximate calculation of $\bW$ can be made to have the effect of 
replacing $\bW$ with a 
matrix $\hat{\bW}\geq \bW$ (i.e., so that $\hat{\bW} - \bW$ is positive 
definite).  To see this, observe that with the replacement of 
$\hat{\bW}$ for $\bW$ (\ref{eq:kahan3}) becomes
\begin{equation}
(\bH-\rho\bI)\hat{\by}=\hat{L}\ \left[(\bH-\rho\bI)- 
\rho(\bI-\rho(\bH^{-1}+\bH^{-1}\bC^{t}\hat{\bW}\bC\bH^{-1}))\right]\hat{\by}
\label{eq:goer3}
\end{equation}
The right-hand side of (\ref{eq:kahan3}) has been replaced with a 
larger right-hand side in (\ref{eq:goer3}).  The left hand side 
remains the same, so (\ref{eq:goer3}) and (\ref{eq:kahan3}) will have the 
same numbers of positive ($\pi$) and negative ($\nu$) eigenvalues.  
The min-max characterization then may be used to show that
\begin{eqnarray*}
L_{k} \leq \hat{L}_{k}< 0 &\mbox{for} & k=1,\ 2,\ \dots,\ \nu \\
0< \hat{L}_{-\ell} \leq L_{-\ell} &\mbox{for} & \ell=1,\ 2,\ \dots,\ \pi.
\end{eqnarray*}
The inequalities of (\ref{LLbnds}) remain valid if $\hat{L}_{k}$ 
replaces $L_{k}$ and $\hat{L}_{-\ell}$ replaces $L_{-\ell}$.
Likewise if we define $\hat{\Lambda}_{\pm i}^{(L)}=-\rho 
\hat{L}_{\mp i}/(1-\hat{L}_{\mp i})$, the usual labeling is retained
$$
 \dots \hat{\Lambda}^{(L)}_{-3}\leq \hat{\Lambda}^{(L)}_{-2} \leq \hat{\Lambda}^{(L)}_{-1} < \rho 
 < \hat{\Lambda}^{(L)}_{1} \leq \hat{\Lambda}^{(L)}_{2} \leq 
 \hat{\Lambda}^{(L)}_{3}\dots,
$$ 
and $\hat{\Lambda}^{(L)}_{-k}\leq \Lambda^{(L)}_{-k}$ for each $k=1,\, 
\dots,\, \nu$.  The situation regarding the positively indexed 
$\hat{\Lambda}^{(L)}$ that yield bounds above $\rho$ is 
slightly more complicated since it may occur that 
$\hat{L}_{-\ell} < 1 < L_{-\ell}$ which would then imply that 
$\hat{\Lambda}_{\ell}^{(L)}<0$.  In effect, 
$\hat{\Lambda}_{\ell}^{(L)}$ has ``wrapped around'' the point at 
infinity, yielding only trivial bounds for $\lambda_{r+\ell-1}$.
Nontrivial bounds are retained whenever 
$\hat{\Lambda}_{\ell}^{(L)}>0$, however.

Now,  much the same development that yielded Theorem \ref{leftdefKahan} may 
be followed with $\hat{\bW}$ replacing $\bW$.  This is summarized as
\begin{theorem} \label{GoerKahan}
Suppose $\rho$ is not an eigenvalue of (\ref{eigprob}).  Each 
interval $[\hat{\Lambda}^{(L)}_{-i},\ \rho)$ and $(\rho,\ 
\hat{\Lambda}^{(L)}_{j}]$ contains respectively at least  $i$ and $j$ 
eigenvalues of (\ref{eigprob}), where
	$$
 0<\hat{\Lambda}^{(L)}_{-\nu} \leq \dots \leq \hat{\Lambda}^{(L)}_{-2} \leq 
 \hat{\Lambda}^{(L)}_{-1}< \rho < \hat{\Lambda}^{(L)}_1 \leq \hat{\Lambda}^{(L)}_2 \leq 
 \dots
	$$
are the positive eigenvalues of the $(m+k) \times (m+k)$ matrix pencil
\begin{align}
	\left[
	\begin{array}{cc}
		\bH & \bC^{t}  \\
		\bC & \hat{\bN}_{1}
	\end{array}\right]& - \Lambda^{(L)}\left[
	\begin{array}{cc}
		\bI & 0  \\
0 & \hat{\bM}_{1} \end{array}\right], \label{eq:goer5} \\
 \mbox{where }\quad \hat{\bM}_{1}& =\frac{1}{\rho}(\hat{\bN}_{1}-\bN_{2}) 
 \nonumber \\
 \hat{\bN}_{1}& =\hat{\bW}^{-1}+\bC\bH^{-1}\bC^{t}, \nonumber \\
	\bN_{2}& =\bC(\bH-\rho\bI)^{-1}\bC^{t}, \nonumber
\end{align}
and $\hat{\bW}$ is any positive-definite matrix satisfying $\hat{\bW}\geq 
\bW=\bQ_{2}^{t}\bK^{-1}\bQ_{2}$.\newline 
$\rho$ is an eigenvalue of (\ref{eq:goer5}) with multiplicity $k$.
\end{theorem}

Goerisch (\cite{Goerhab},\cite{GoeHau85},\cite{GoeZim86}) discovered 
this approach and developed a very flexible 
framework for applying this critical approximation step for the 
original left-definite Lehmann formulation (\ref{LdefPform}) in a PDE 
setting.  He called 
it the $\{ \mathcal{X},\, b, T\}$ method (referring to an auxiliary 
vector space $\mathcal{X}$, an auxiliary bilinear form $b$, and an 
auxiliary linear operator $T$ that he introduces) but most others 
refer to this approach simply as the Lehmann-Goerisch method.  To give 
a simple example, suppose a lower bound to $\bK$ is known: $\kappa 
\|\bx\|^{2}\leq \bx^{t}\bK\bx$, and suppose we have obtained an 
approximate solution $\bZ_{2}$ to the matrix equation $\bK\bZ=\bQ_{2}$.  Let 
$\bR=\bQ_{2}-\bK\bZ_{2}$ be the associated residual matrix.  Then one 
may verify that
\begin{align*}
\bW=\bQ_{2}^{t}\bK^{-1}\bQ_{2}& =\bR^{t}\bK^{-1}\bR 
+\bZ_{2}^{t}\bR+\bQ_{2}^{t}\bZ_{2} \\
	 & \leq \frac{1}{\kappa}\bR^{t}\bR 
	 +\bZ_{2}^{t}\bR+\bQ_{2}^{t}\bZ_{2}\stackrel{def}{=}\hat{\bW}.
\end{align*}
Note that $\hat{\bW}$ contains the nominal estimate of $\bW$, 
$\bQ_{2}^{t}\mbox{``}\bK^{-1}\bQ_{2}\mbox{''}=\bQ_{2}^{t}\bZ_{2}$,
together with correction terms that ensure $\hat{\bW}\geq \bW$ and
that can be made small by solving 
$\bK\bZ=\bQ_{2}$ more accurately.

\section{A Left-Right Comparison}
For the general eigenvalue problem (\ref{eigprob}), application of 
either right- or left-definite Lehmann bounds involve solving linear 
systems having either $\bM$ (for right-definite problems) or $\bK$ (for 
left-definite problems) as a coefficient matrix.  If one system is very much 
simpler than the other (e.g., if $\bM=\bI$) one may feel compelled to 
choose the simpler path.  But is there a difference in accuracy ?  
Goerisch and coworkers in Braunschweig and 
Clausthal (see for example, \cite{GoeAlb83} and \cite{GoeHau85})
have observed  that for many applications in PDE settings, 
left-definite Lehmann bounds often were superior to right-definite 
bounds --- even if an extra level of approximation is included as
 described in Theorem \ref{GoerKahan}. Along similar lines,
Knyazev \cite{Kny} has produced error estimates for Lehmann methods that suggest 
left-definite bounds might be better than right-definite bounds 
asymptotically.

We explore this issue here.   Define
$$
\bJ_{0}=\bH_{0}-\rho \bH_{1},\ \bJ_{1}=\bH_{1}-\rho \bH_{2},\ 
\mbox{and }\ \bJ_{2}=\bH_{2}-\rho \bH_{3}.
$$
The matrix pencils associated with (\ref{RdefPform}) and (\ref{LdefPform}) may be 
written as
\begin{equation}
	\bJ_{1}-R(\bJ_{0}-\rho \bJ_{1})
	\label{RdefPformshort}
\end{equation}
and 
\begin{equation}
	\bJ_{1}-L(\bJ_{1}-\rho \bJ_{2})
	\label{LdefPformshort}
\end{equation}
for right-definite and left-definite problems, respectively.

The following lemma and theorem incorporate some unpublished results of 
Goerisch\footnote{Friedrich Goerisch died suddenly in 1995 
after a brief illness. The loss of his passion and insight is still 
deeply felt among his colleagues and friends.}.
\begin{lemma}
	Let $\bG= \left[\begin{array}{cc}
		\bJ_{0} & \bJ_{1}  \\
		\bJ_{1} & \bJ_{2}
	\end{array}\right]\in {\mathbb {R}}^{2m \times 2m}$. $\bG$ has no more 
	than $r-1$ negative eigenvalues.
\end{lemma}
\textsc{ Proof:} Suppose that $\bG$ has $r$ or more negative eigenvalues. 
Then there is an $r$-dimensional subspace $\mathcal{Z}$ of  
${\mathbb {R}}^{2m}$ such that $\bz^{t}\bG\bz<0$ for all $\bz\in \mathcal{Z}$
with $\bz\neq 0$.  Define the linear mapping 
$T:\mathcal{Z}\rightarrow {\mathbb {R}}^{n}$ by 
$$
T(\bz)=\sum_{i=1}^{m} z_{i}\bK\bp_{i}+\sum_{i=1}^{m}
z_{i+m}\bM\bp_{i}
$$
Elementary manipulations verify that for  $\bz\in \mathcal{Z}$
with $\bz\neq 0$,
\begin{equation}
	\bz^{t}\bG\bz=T(\bz)^{t}\bM^{-1}T(\bz)-\rho T(\bz)^{t}\bK^{-1}T(\bz)<0
	\label{lem_ass1}
\end{equation}
In particular, this means that $T(\bz)=0$ implies that $\bz=0$, so 
$null(T)=0$ and $rank(T)=\dim \mathcal{Z} =r$.  

Since $\bK$ is positive-definite $\bu^{t}\bK^{-1}\bu >0$ for all $\bu 
\in {\mathbb {R}}^{n}$ so (\ref{lem_ass1}) implies 
$\bu^{t}\bM^{-1}\bu/\bu^{t}\bK^{-1}\bu <\rho$ for all $\bu\in Ran(T)$ 
with $\bu\neq0$.  

Now $\lambda$ is an eigenvalue of (\ref{eigprob}) if and only if it is also an 
eigenvalue of $\bM^{-1}\bv=\lambda \bK^{-1}\bv$, so by the min-max 
principle
$$
\lambda_{r}=\min_{\dim \mathcal{P}=r} \max_{\bu\in\mathcal{P}}
\frac{ \bu^{t}\bM^{-1}\bu}{\bu^{t}\bK^{-1}\bu } \leq 
\max_{\bu\in Ran(T)}
\frac{ \bu^{t}\bM^{-1}\bu}{\bu^{t}\bK^{-1}\bu } < \rho
$$
which contradicts $\lambda_{r-1}<\rho<\lambda_{r}$.  Thus, $\dim 
\mathcal{Z} < r$.    $\quad \blacksquare$

\begin{theorem} 
	If the harmonic 
Ritz value $\tilde{\Lambda}_{r-1}$ from (\ref{HarmRitzeig}) satisfies 
$\tilde{\Lambda}_{r-1}<\rho$ then left-definite Lehmann bounds will be 
uniformly better than right-definite Lehmann bounds:
 \begin{equation}
	 \Lambda^{(R)}_{-k} \leq \Lambda^{(L)}_{-k}\leq \lambda_{r-k}
	 \quad \mbox{ for }k=1,\ \dots,\ r-1
 	\label{LMconclusion1}
 \end{equation}
\begin{equation}
	\lambda_{r+\ell-1} \leq  \Lambda^{(L)}_{\ell} \leq \Lambda^{(R)}_{\ell}
	\quad \mbox{ for }\ell=1,\ \dots,\ m-r+1
	\label{LMconclusion2}
\end{equation}
\end{theorem}

 \textsc{Proof}:
To show that (\ref{LMconclusion1}) and (\ref{LMconclusion2}) are true, it is sufficient to show 
that  $L_{k}\leq 1+\rho R_{k}$ for $k=1,\ 2,\ \dots, r-1$ and that 
$1+\rho R_{-\ell}\leq L_{-\ell}$ for $\ell=1,\ 2,\ \dots, m-r+1$
From (\ref{RdefPformshort}), one finds that $1+\rho R_{k}$ and $1+\rho 
R_{-\ell}$ are  eigenvalues of 
\begin{equation}
	\bJ_{0}-(1+\rho R)(\bJ_{0}-\rho \bJ_{1}).  
	\label{eq:goeprf1}
\end{equation}
Since $\Lambda_{r-1}\leq \tilde{\Lambda}_{r-1}<\rho$, both $\bJ_{0}$ 
and  $\bJ_{1}$ have $r-1$ negative eigenvalues.  This 
implies that both 
(\ref{RdefPformshort}) and  
(\ref{LdefPformshort}) have $r-1$ negative eigenvalues.
Premultiplication of (\ref{eq:goeprf1}) by $\bJ_{1}\bJ_{0}^{-1}$ yields an equivalent 
matrix pencil:
$$
\bJ_{1}-(1+\rho R)(\bJ_{1}-\rho \bJ_{1}\bJ_{0}^{-1}\bJ_{1})
$$

Consider 
$$\bG= \left[\begin{array}{cc}
		\bJ_{0} & \bJ_{1}  \\
		\bJ_{1} & \bJ_{2}
	\end{array}\right]=\left[\begin{array}{cc}
		\bI & 0  \\
		\bJ_{1}\bJ_{0}^{-1} & \bI
	\end{array}\right]\left[\begin{array}{cc}
		\bJ_{0} & 0 \\
		0 & \bJ_{2}-\bJ_{1}\bJ_{0}^{-1}\bJ_{1}
	\end{array}\right]\left[\begin{array}{cc}
		\bI & \bJ_{0}^{-1}\bJ_{1}  \\
		0 & \bI
	\end{array}\right].
	$$
By the lemma and the Sylvester law of inertia, 
$\bJ_{0}\oplus \bJ_{2}-\bJ_{1}\bJ_{0}^{-1}\bJ_{1} $ 
can have no more than $r-1$ negative eigenvalues. Since 
$\bJ_{0}$ has exactly $r-1$ eigenvalues by hypothesis,
$\bJ_{2}-\bJ_{1}\bJ_{0}^{-1}\bJ_{1}$ must be positive semi-definite
and 
$$
0< \bx^{t}(\bJ_{1}-\rho \bJ_{2})\bx
<\bx^{t}(\bJ_{1}-\rho \bJ_{1}\bJ_{0}^{-1}\bJ_{1})\bx
$$ 
for all nontrivial $\bx$.
Hence, for $k=1,\ 2,\ \dots,\ r-1$,
\begin{align*}
	1+\rho R_{k}=& \min_{\dim {\mathcal S} = k}\max_{\bx\in {\mathcal S}}
	\frac{\bx^{t}\bJ_{1}\bx}{\bx^{t}(\bJ_{1}-\rho 
	\bJ_{1}\bJ_{0}^{-1}\bJ_{1})\bx} \\
	\geq & \min_{\dim {\mathcal S} = k}\max_{\bx\in {\mathcal S}}
	\frac{\bx^{t}\bJ_{1}\bx}{\bx^{t}(\bJ_{1}-\rho \bJ_{2})\bx} =L_{k}
\end{align*}
and for $\ell=1,\ 2,\ \dots,\ m-r+1$,
\begin{align*}
	-(1+\rho R_{-\ell})= & \min_{\dim {\mathcal S} = \ell}\max_{\bx\in {\mathcal S}}
	\frac{-\bx^{t}\bJ_{1}\bx}{\bx^{t}(\bJ_{1}-\rho 
	\bJ_{1}\bJ_{0}^{-1}\bJ_{1})\bx} \\
	\geq & \min_{\dim {\mathcal S} = \ell}\max_{\bx\in {\mathcal S}}
	\frac{-\bx^{t}\bJ_{1}\bx}{\bx^{t}(\bJ_{1}-\rho \bJ_{2})\bx} =-L_{-\ell}.
\end{align*}
Since there will be subspaces of dimension up to $r-1$ for 
which  $\bx^{t}\bJ_{1}\bx<0$ and subspaces of dimension up to $m-r+1$ for 
which  $\bx^{t}\bJ_{1}\bx>0$,  we may restrict ourselves  to $\bx$ for which 
the numerators in the above expressions are strictly negative with 
no loss of generality.  $ \quad \blacksquare$

\section{A Ritz-Lehmann Comparison}
One may hope that the role spectral mapping played in the derivation
of both left- and right-definite variants of Lehmann's method might 
lead to significant improvements beyond the straightforward application 
of the Rayleigh-Ritz method. Indeed, spectral mapping has been used for 
some time with Lanczos methods (e.g., \cite{ErRu80}) with sometimes 
spectacular effect and so encouraged, some have 
considered the use of right-definite Lehmann bounds using Krylov 
subspaces generated in the course of an ordinary Lanczos process (e.g., 
\cite{Mor} and \cite{PaiParVor}).    By and large, results along these 
lines have been disappointing when compared with what 
``shift-and-invert'' methods offer (albeit at a much higher price). 
One may instead seek to  compare the expected outcomes of 
Lehmann methods with those of Rayleigh-Ritz methods.  Observe 
that each method makes optimal use of the information 
required in the sense that no better bounds are possible with the 
infomation used, so in a certain manner of speaking we are really comparing 
the utility of various types of information in extracting eigenvalue 
information.

Zimmerman \cite{Zimm}  proved that the error in left-definite Lehmann bounds is no 
worse than proportional to the error in Ritz bounds and may be 
smaller.  Thus, left-definite Lehmann bounds carry the potential of greater 
accuracy than Ritz bounds. We probably shouldn't expect them to be 
 much better, though. In  \cite{Kny}, Knyazev states that eigenvector approximations 
provided by either the right- or left-definite variants of Lehmann's 
method will asymptotically approach 
the corresponding Ritz vectors as they close upon the true 
eigenvectors.  Thus, Lehmann methods appear to recover invariant 
subspace information with about the same efficiency as Rayleigh-Ritz 
methods.

It is important to note that Lehmann methods provide eigenvalue 
\emph{bounds} that often are difficult to obtain in other ways.
For example, Behnke  
\cite{Beh88} combined right-definite Lehmann methods with interval techniques
 in order to deduce guaranteed bounds to matrix eigenvalue 
problems and his approach appears to be competitive with the best known 
interval algorithms for this problem.

For the remainder of this section, we will consider 
the application of a left-definite Lehmann method within a 
Lanczos process for resolving a large-scale matrix eigenvalue problem. 
Since left-definite Lehmann methods are known to be superior to 
right-definite Lehmann methods (at least to the extent claimed in 
Section 4), one may seek to improve upon the results of Morgan 
\cite{Mor} by using left-definite Lehmann-Goerisch bounds as formulated 
in Theorem \ref{GoerKahan}. 

Specifically, let $\bM=\bI$ in (\ref{eigprob}) and let $\bT$ be a 
tridiagonal matrix that is similar to $\bK$ -- so that $\bK=\bQ \bT 
\bQ^{t}$ for some $n \times n$ unitary matrix $\bQ$.  For any index $1 
\leq \ell \leq n $, let $\bT_{\ell}$ denote the $\ell$th principal 
submatrix of $\bT$:
	$$
     \bT_{\ell}=\left[	\begin{array}{cccccc}
		\alpha_{1} & \beta_{1}&  &  &  &   \\
		\beta_{1} & \alpha_{2} & \beta_{2} &  &  &   \\
		 & \beta_{2} & \alpha_{3} & \ddots &  &   \\
		 &  & \ddots & \ddots &  &   \\
		 &  &  &  &  & \beta_{\ell-1} \\
		 &  &  &  & \beta_{\ell-1} & \alpha_{\ell}
	\end{array}  \right]
	$$
 and define $\bV$ via a partitioning of $\bT$ as
	$$
	\bT=\left[ 
	\begin{array}{cc}
		\bT_{\ell} & \beta_{\ell}\be_{\ell}\be_{1}^{t}  \\
		\beta_{\ell}\be_{1}\be_{\ell}^{t}  &\bV
	\end{array}\right]
	$$
	Let $\bQ_{\ell}$ denote a matrix containing the first $\ell$ columns 
	of $\bQ$: $\bQ_{\ell}=\left[ \bq_{1},\ \ldots \ \bq_{\ell} \right]$.
	
	The Lanczos algorithm builds up the 
	matrices $\bT$ and $\bQ$  one column at a time starting with 
	the vector $\bq_{1}$. Only information on the 
	action of $\bK$  on selected vectors in $\mathbb{R}^{n}$ 
	is used.  	Different choices for $\bq_{1}$ 
	 produce distinct outcomes for $\bT$, if all goes 
	well.  Extracting useful information when not all goes well
	is fundamental to modern approaches -- a discussion may be found in 
	\cite{Parlett}.
	
	At the $\ell$th step, the basic Lanczos recursion appears as
		$$ \bK \bQ_{\ell}= 
		\bQ_{\ell}\bT_{\ell}+\beta_{\ell}\bq_{\ell+1}\be_{\ell}^{t} $$
 In exact arithmetic, the first $\ell$ 
steps yields a matrix $\bQ_{\ell}$ that satisfies
$\bQ_{\ell}^{t}\bQ_{\ell}=\bI$ and
	$$
	Ran(\bQ_{\ell})=  \mbox{span}\{\bq_{1},\, \bA\bq_{1},\, \dots,\, 
	\bA^{\ell-1}\bq_{1}\}={\mathcal{K}}_{\ell}(\bA,\bq_{1}),
	$$
a Krylov subspace of order $\ell$. The application of Theorem 
\ref{GoerKahan} is straightforward: 
\begin{theorem} Let $\bM=\bI$ and suppose $\rho$ 
is not an eigenvalue of (\ref{eigprob}).  Each interval 
$[\Lambda^{(L)}_{-i},\ \rho)$ and $(\rho,\ \Lambda^{(L)}_{j}]$ 
contains respectively at least $i$ and $j$ eigenvalues of the matrix $\bK$, where
	$$
 0< \Lambda^{(L)}_{-\nu}\dots \leq \Lambda^{(L)}_{-2} \leq 
 \Lambda^{(L)}_{-1}< \rho< \Lambda^{(L)}_1 \leq \Lambda^{(L)}_2 \leq 
 \dots
	$$
are the positive eigenvalues of the tridiagonal matrix pencil
	\begin{equation}
		\left[\begin{array}{cc}
			 \bT_{\ell} & \beta_{k}\be_{k}  \\
	         \beta_{k}\be_{k}^{t} & \omega_{k+1}^{-1} + \beta^{2}_{k}\be_{k}^{t} 
	\bT_{k}^{-1}\be_{k} \end{array}\right]-\Lambda^{(L)}\left[\begin{array}{cc} 
	\bI & 0 \\
	         0 & (\rho\omega_{k+1})^{-1} - \beta^{2}_{k}\delta_{k+1}(\rho) \end{array}\right]
		\label{LanGoerKah}
	\end{equation}
	where $\omega_{k+1}$ is any number that satisfies
	\begin{align*}
\omega_{k+1}\geq \bq_{k+1}\bK^{-1}\bq_{k+1}\\
			\mbox{and} \quad \delta_{k+1}(\rho) & =\be_{k}^{t} 
			\bT_{k}^{-1}(\bT_{k}-\rho)^{-1}\be_{k}.
		\end{align*}
Note that $\rho$ is a simple eigenvalue of (\ref{LanGoerKah})
\end{theorem}

We apply this directly to the numerical example considered in \cite{PaiParVor} 
and in Section 1.  Figure 2 shows the convergence history both for  
Ritz bounds and for left-definite Lehmann bounds, for the seventh 
through tenth eigenvalues of the matrix. We also apply a shift and 
invert Lanczos method using the spectral transformation $\lambda \mapsto 
\frac{\lambda}{\lambda - \rho}$.  A few features are apparent.
The first is that the Lehmann bounds aren't nearly as good as the 
shift and invert bounds to which they are closely related.  Paige, 
Parlett, and van der Vorst \cite{PaiParVor} observed this 
disappointing behaviour for right-definite Lehmann methods (in their 
context, harmonic Ritz on a shifted matrix) --- the left-definite Lehmann 
method does not fare much better.  Knyazev's observations \cite{Kny} relating  
convergence of Lehmann eigenvectors to Ritz vectors suggest that 
spectral information for interior matrix eigenvalues will not be 
picked up any more rapidly with Lehmann methods than for Ritz methods.  
This is in stark contrast with shift and invert strategies which will 
produce approximate eigenvectors that are rapidly drawn into invariant subspaces 
associated with eigenvalues close to $\rho$.  

The second observation 
is that, nonetheless, the Lehmann bounds do appear to approach the 
exact eigenvalues at a rate comparable to that of the Ritz bounds --- 
consistent with the results of Zimmerman discussed above.  Furthermore,
one can see that 
the Lehmann bounds appear to pass through a series of stagnation 
points en route to their limit, and the farther they lie from $\rho$, the more 
abrupt the transition between stagnation points.  These 
stagnation points appear to be close to the exact matrix eigenvalues.

\begin{figure}
  \centerline{\psfig{figure=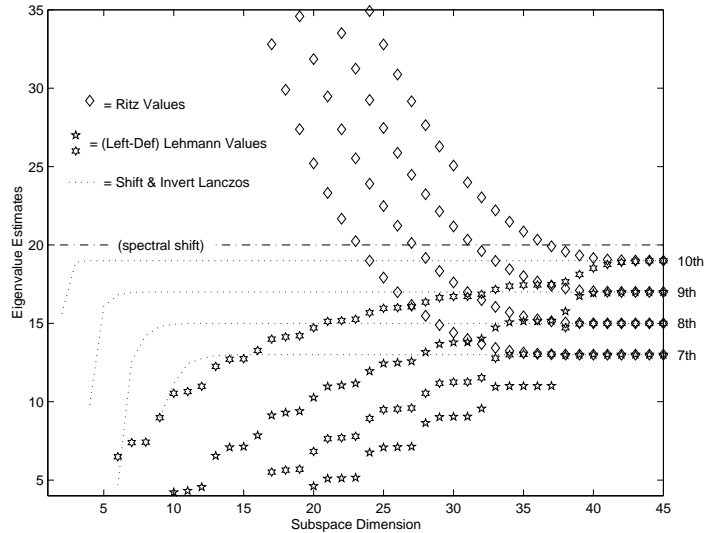,height=3in}}
 \caption{ Convergence of Ritz and Lehmann bounds using 
		Krylov subspaces vs. Shift \& invert Lanczos with same 
		starting vector.}
 \label{fig:Ritz3}
 \end{figure}
 
The following simple Bauer-Fike style perturbation result lends some 
insight to this behaviour.

\begin{theorem}
	Let $\Lambda^{(L)}$ be any left-definite Lehmann bound and denote with
	$\Lambda_{i}$ the Ritz values from (\ref{Ritzeig}).  Then
	\begin{equation}
		\min_{i}\left(\frac{|\Lambda_{i}-\rho|}{\rho}\right)\left(\frac{|\Lambda_{i}-\Lambda^{(L)}|}{\Lambda^{(L)}}\right)
		\Lambda_{i}\quad \leq \quad \|\bW\|\, \|\bC\|^{2}
\label{eq:bauerfike}
	\end{equation}

\end{theorem}
\textsc{Proof:} If either $(\bH-\rho\bI)$ or $(\bH-\Lambda^{(L)}\bI)$ 
is singular then (\ref{eq:bauerfike}) holds trivially.  Suppose 
then that $(\bH-\rho\bI)$ and $(\bH-\Lambda^{(L)}\bI)$ are 
nonsingular.  Rearrange the expression (\ref{eq:kahan4}) to get
$$
\by=-\rho 
\Lambda^{(L)}(\bH-\rho\bI)^{-1}(\bH-\Lambda^{(L)}\bI)^{-1}\bH^{-1}\bC^{t}\bW\bC\by
$$
Take norms on each side and simplify:
 \begin{equation}
	 1\leq\rho 
	\Lambda^{(L)}\|(\bH-\rho\bI)^{-1}(\bH-\Lambda^{(L)}\bI)^{-1}\bH^{-1}\|\; 
	\|\bW\|\,\|\bC\|^{2}
 \label{eq:BF1}
 \end{equation}
Then notice that
\begin{align*}
	\|(\bH-\rho\bI)^{-1}(\bH-\Lambda^{(L)}\bI)^{-1}\bH^{-1}\|= & \max_{i}
	\left(\frac{1}{|\Lambda_{i}-\rho|}\right)\left(\frac{1}{|\Lambda_{i}-\Lambda^{(L)}|}\right)
	\frac{1}{\Lambda_{i}} \\
= & 1/\min_{i} (|\Lambda_{i}-\rho|\, |\Lambda_{i}-\Lambda^{(L)}|\, \Lambda_{i}),
\end{align*}
which may be combined with (\ref{eq:BF1}) to get (\ref{eq:bauerfike}).
$\quad \blacksquare$

Notice that the right hand side of (\ref{eq:bauerfike}) has a magnitude 
related to the size of the Ritz residual  $\bK\bQ_{1}-\bQ_{1}\bH$ 
and is independent of which Lehmann bound $\Lambda^{(L)}$ is chosen.  
Suppose the right hand side of (\ref{eq:bauerfike}) is moderately 
small and choose a Lehmann bound $\Lambda^{(L)}$.  If $\Lambda^{(L)}$ 
is not close to $\rho$ then any Ritz value $\Lambda_{i}$ that is close to 
$\Lambda^{(L)}$ will not be close to $\rho$ either.  Thus any $\Lambda^{(L)}$
chosen far from $\rho$ is 
constrained by (\ref{eq:bauerfike}) to be nearer to 
 at least one $\Lambda_{i}$ then it would be were $\Lambda^{(L)}$ chosen closer to 
$\rho$.  A qualitative interpretation that one might take from this is 
that Lehmann bounds $\Lambda^{(L)}$ far from $\rho$ tend to occur in 
the neighborhood of Ritz values $\Lambda_{i}$.  Furthermore, Lehmann 
bounds $\Lambda^{(L)}$ far from $\rho$ that are also situated toward 
the edges of the spectrum will tend to aggregate in the neighborhood 
of exact eigenvalues since the attracting Ritz values themselves will 
be approximating extreme eigenvalues fairly well.

\end{document}